\documentclass{amsproc}
\usepackage{epsfig,amscd,amssymb,amsmath,amsfonts}
\usepackage[margin=1.2in]{geometry}
\newtheorem{theorem}{Theorem}[section]

\theoremstyle{definition}
\newtheorem{definition}[theorem]{Definition}

\newtheorem{example}[theorem]{Example}

\theoremstyle{remark}
\newtheorem{remark}[theorem]{Remark}

\numberwithin{equation}{section}

\begin{document}

\title{Higher Order and Secondary Hochschild Cohomology}

\author{Bruce R. Corrigan-Salter}
\address{Department of Mathematics,
Wayne State University,
Detroit, MI 48202, USA}
\email{brcs@wayne.edu}

\author{Mihai D. Staic}
\address{Department of Mathematics and Statistics, Bowling Green State University, Bowling Green, OH 43403 }
\email{mstaic@bgsu.edu}



\subjclass[2010]{Primary  16E40, Secondary   18G30}
\date{January 1, 1994 and, in revised form, June 22, 1994.}


\keywords{Hochschild Cohomology}

\begin{abstract} In this note we give a generalization for the higher order Hochschild cohomology and show that the secondary Hochschild cohomology is a particular case of this new construction.  
\end{abstract}

\maketitle


%

\section*{Introduction}

Hochschild cohomology is a useful tool for studying deformation theory and  it was studied extensively over the years (for example see \cite{MR0161898}, \cite{MR0171807}, \cite{MR981619}, \cite{MR551624}, \cite{MR1217970} and \cite{MR2153227}).

Higher order Hochschild (co)homology was introduced by Pirashvili in \cite{MR1755114} (see also \cite{MR0339132} and  \cite{MR2383113}). It associated to a simplicial set $X_{\bullet}$, a commutative  $k$-algebra $A$ and an $A$-bimodule $M$, the higher Hochschild cohomology groups $H^n_{X_{\bullet}}(A,M)$. When $X_{\bullet}$ is the standard simplicial set associated to the sphere $S^1$, one recovers the usual Hochschild cohomology. One  important feature of these cohomology groups is that they depend only on the homotopy type of the geometric  realization of the simplicial set $X_{\bullet}$. For more recent results about higher Hochschild cohomology see \cite{MR3338542} and \cite{MR3173402}.

Secondary Hochschild cohomology was introduced in \cite{MR3465889} where it was used to study $B$-algebra structures on the algebra $A[[t]]$. It associates to a triple $(A,B,\varepsilon)$ (where $\varepsilon$ gives the $B$-algebra structure on $A$), and an $A$-bimodule $M$ that is $B$-symmetric, the secondary Hochschild cohomology groups $H^n((A,B,\varepsilon),M)$.  The main difference from the usual Hochschild cohomology is that the secondary Hochschild cohomology gives information about both, the product and the $B$-algebra structure on $A[[t]]$. Many of the results that are known for Hochschild cohomology can be extended to secondary Hochschild cohomology. One can study $B$-extensions of the algebra $A$, show that $H^{\bullet}((A,B,\varepsilon), A)$ is a multiplicative operad and admits a Hodge type decomposition (\cite{MR3338544}). There is also a cyclic version of the secondary cohomology that can be described using a generalization of the bar resolution (\cite{s4}).

Our main goal in this paper is to show that secondary Hochschild cohomology is a certain version of higher order Hochschild cohomology.  More precisely, we consider a simplicial pair $(X_{\bullet},Y_{\bullet})$ (where $Y_{\bullet}$ is a simplicial set and $X_{\bullet}$ is simplicial subset of $Y_{\bullet}$), a triple $(A,B,\varepsilon)$ (where $A$ and $B$ are  commutative $k$-algebras,  and $\varepsilon:B\to A$ is a morphism of $k$-algebras) and $M$ a symmetric  $A$-bimodule. To this setting we associate the groups $H^n_{(X_{\bullet},Y_{\bullet})}((A,B,\varepsilon),M)$. When $X_{\bullet}=Y_{\bullet}$ we recover the  higher order Hochschild cohomology $H^n_{X_\bullet}(A,M)$. When $(X_{\bullet},Y_{\bullet})=(S^1, D^2)$ with the natural simplicial structure, we recover the secondary Hochschild cohomology $H^n((A,B,\varepsilon), M)$.

\section{Preliminaries}

In this paper we fix a field $k$ and denote $\otimes_k$ by $\otimes$.  We assume that the reader is familiar with Hochschild cohomology, but provide some details for the discussion of higher order and secondary Hochschild cohomology.  We also assume familiarity with simplicial sets.  

\subsection{Higher order Hochschild cohomology}
 We follow the description in \cite{MR2383113} (see also  \cite{MR1755114}).  Assume that $A$ is a commutative $k$-algebra and $M$ is a symmetric $A$-bimodule.

Let $V$ be a finite pointed set such that $\vert V\vert =v+1$ (we  identify it with $v_+=\{ 0,1,...,v\}$  with $0$ the fixed element) and define $\mathcal{H}(A,M)(V)=\mathcal{H}(A,M)(v_+)=Hom_k(A^{\otimes v},M)$. 
For $\phi:V=v_+\to W=w_+$   we define 
$$\mathcal{H}(A,M)(\phi): \mathcal{H}(A,M)(w_+)\to \mathcal{H}(A,M)(v_+)$$
determined as follows: if  $f\in \mathcal{H}(A,M)(w_+)$
then $$\mathcal{H}(A,M)(\phi)(f)(a_1\otimes ...\otimes a_v)=b_0f(b_1\otimes ...\otimes b_{w})$$ 
where $$b_i=\prod_{\{j\in V | j\neq 0, \phi (j)=i\}}a_j.$$

Take $X_{\bullet}$ to be a pointed simplicial set. Suppose that $\vert X_n\vert =s_n+1$, we identify the set $X_n$ with $(s_{n})_{+}=\{ 0,1,...,s_n\}$ then define 
$$C_{X_\bullet}^n=\mathcal{H}(A,M)(X_n)=Hom_k(A^{\otimes s_n},M).$$ For each $d_i:X_{n+1}\to X_{n}$ we define $d_i^*=\mathcal{H}(A,M)(d_i):C_{X_{\bullet}}^{n}\to C_{X_{\bullet}}^{n+1}$ and take $\partial_{n}:C^{n}_{X_{\bullet}}\to C^{n+1}_{X_{\bullet}}$ defined as 
$\partial_n=\sum_{i=0}^{n+1}(-1)^i(d_i)^*$. 

The homology of this complex is denoted by $H^n_{X_{\bullet}}(A,M)$ and is called the higher order Hochschild cohomology group. One  interesting fact is that these groups depend only on the homotopy type of the geometric realization of the simplicial set $X_{\bullet}$. 

\subsection{Secondary Hochschild cohomology}
We recall the construction from \cite{MR3465889}. Let $A$ be a $k$-algebra, $B$ a commutative $k$-algebra, $\varepsilon:B\to A$ a morphism of $k$-algebras such that $\varepsilon(B)\subset \mathcal{Z}(A)$  and $M$ an $A$-bimodule that is $B$-symmetric. 

We  define $C^n((A,B,\varepsilon);M):=Hom_k(A^{\otimes n}\otimes B^{\otimes \frac{n(n-1)}{2}},M)$ and $$\delta_{n-1}^{\varepsilon}:C^{n-1}((A,B,\varepsilon);M)\to C^{n}((A,B,\varepsilon);M)$$ such that for $f\in C^{n-1}((A,B,\varepsilon);M)$, we have:
\begin{equation}
\delta^{\varepsilon}_{n-1}(f)\left(\displaystyle\otimes
\left(
\begin{array}{cccccccc}
 a_{1}& \alpha_{1,2} &  ...&\alpha_{1,n-1}&\alpha_{1,n}\\
1 & a_{2}       & ...&\alpha_{2,n-1}&\alpha_{2,n}\\
. & .    &...&.&.\\
1& 1& ...&a_{n-1}&\alpha_{n-1,n}\\
1 & 1&...&1&a_{n}\\
\end{array}
\right)\right)=\label{delta}
\end{equation}
\begin{eqnarray*}
&
a_1\varepsilon(\alpha_{1,2}\alpha_{1,3}...\alpha_{1,n})f\left(\displaystyle\otimes
\left(
\begin{array}{cccccc}
 a_{2}       & \alpha_{2,3} &...&\alpha_{2,n-1}&\alpha_{2,n}\\
1  & a_{3} &...    &\alpha_{3,n-1}&\alpha_{3,n}\\
 .       &. &...&.&.\\
1& 1  &...&a_{n-1}&\alpha_{n-1,n}\\
 1& 1  &...&1&a_{n}\\
\end{array}
\right)\right)+&\\
&\;&\;\\
& \sum_{i=1}^{n-1}(-1)^if\left(\displaystyle\otimes
\left(
\begin{array}{cccccccc}
 a_{1}       & \alpha_{1,2}&... &\alpha_{1,i}\alpha_{1,i+1}&...&\alpha_{1,n-1}&\alpha_{1,n}\\
1  &a_{2}&...&\alpha_{2,i}\alpha_{2,i+1}     &...&\alpha_{2,n-1}&\alpha_{2,n}\\
 .       &. &...&.&...&.&.\\
  1       &1 &...&a_{i}a_{i+1}\varepsilon(\alpha_{i,i+1})&...&\alpha_{i,n-1}\alpha_{i+1,n-1}&\alpha_{i,n}\alpha_{i+1,n}\\
   .       &. &...&.&...&.&.\\
1& 1  &...&.&...&a_{n-1}&\alpha_{n-1,n}\\
 1& 1  &...&.&...&1&a_{n}\\
\end{array}
\right)\right)+&\\
& \;&\\
&(-1)^{n}f\left(\displaystyle\otimes
\left(
\begin{array}{cccccc}
 a_{1}       & \alpha_{1,2} &...&\alpha_{1,n-2}&\alpha_{1,n-1}\\
1  &a_{2} &...     &\alpha_{2,n-2}&\alpha_{2,n-1}\\
 .       &. &...&.&.\\
1& 1  &...&a_{n-2}&\alpha_{n-2,n-1}\\
 1& 1  &...&1&a_{n-1}\\
\end{array}
\right)\right)a_{n}\varepsilon(\alpha_{1,n}\alpha_{2,n}...\alpha_{n-1,n}).&
\end{eqnarray*}
where $a_i\in A$ and $\alpha_{i,j}\in B$. 
It was proved in \cite{MR3465889} that $(C^n((A,B,\varepsilon);M),\delta_n^{\varepsilon})$ is a complex. The homology of this complex is called the secondary Hochschild cohomology and is denoted by $H^n((A,B,\varepsilon);M)$. The homology and the cyclic version of this theory were discussed in \cite{s4}.

\section{Main Construction}

In this section we introduce a new cohomology associated to a simplicial pair $(X_{\bullet},Y_{\bullet})$ and a triple $(A,B,\varepsilon)$ (where $A$ and $B$ are commutative $k$-algebra and $\varepsilon :B\to A$ is a morphism of $k$-algebras).  We start with a few  notations. 

\begin{definition} We consider $\Gamma_{2}$ to be the category whose objects are pairs $(U,V)$, where $V$ is a finite pointed set with base point $*$, and $U$ is a pointed subset of $V$. A morphism $f\in Hom_{\Gamma_{2}}((U_1,V_1), (U_2,V_2))$ is a map of pointed sets $f:V_1\to V_2$ such that $f(U_1)\subset U_2$. 
\end{definition}

\begin{remark} The category of finite pointed sets $\Gamma$ can be see as a full subcategory of $\Gamma_{2}$ in two different ways. First we can take the inclusion given by $V\to (V,V)$, second we can take the inclusion $V\to (\{*\},V)$.
\end{remark}

\begin{definition} A $\Gamma_{2}$-module is a functor from $\Gamma_{2}^{op}$ to $k$-modules.  
\end{definition}
\begin{example}
Let $A$ and $B$ be two commutative $k$-algebras, $\varepsilon:B\to A$ a morphism of $k$-algebras and $M$ a symmetric $A$-bimodule. We construct $${\mathcal L}((A,B,\varepsilon);M):\Gamma_{2}^{op}\to k-mod $$ to be the   $\Gamma_{2}$-module determined as follows.  
For $(U,V)\in \Gamma_{2}$  such that $|U|=1+m$ and $|V|=1+m+n$ define 
$${\mathcal L}((A,B,\varepsilon);M)((U,V))=Hom_k(A^{\otimes m}\otimes B^{\otimes n}, M).$$
If $f:(U_1,V_1)\to (U_2,V_2)$ is a morphism in $\Gamma_2$, we define  $${\mathcal L}((A,B,\varepsilon);M)(f):Hom_k(A^{\otimes m_2}\otimes B^{\otimes n_2}, M) \to Hom_k(A^{\otimes m_1}\otimes B^{\otimes n_1}, M),$$ 
for $\psi\in Hom_k(A^{\otimes m_2}\otimes B^{\otimes n_2}, M)$ then 

${\mathcal L}((A,B,\varepsilon);M)(f)(\psi)(a_1\otimes ...\otimes a_{m_1}\otimes \alpha_1\otimes ...\otimes \alpha_{n_1})=
b_0\cdot \psi(b_1\otimes ...\otimes b_{m_2}\otimes \beta_1\otimes ...\otimes \beta_{n_2})$ where for $i\in U_2$ we have

 \begin{equation} \label{eq} b_i=\prod_{\{j\in U_1|j\neq *,\; f(j)=i\}}a_{j}\prod_{\{k\in V_1\setminus U_1|k\neq *, \; f(k)=i\}}\varepsilon(\alpha_k)\in A,\end{equation}
and for $p\in V_2\setminus U_2$ we have
$$\beta_p=\prod_{\{q\in V_1\setminus U_1|q\neq *, \; f(q)=p\}}\alpha_{q}\in B. $$
With the convention that if the product is taken over the empty set then we put $b_i=1\in A$ and $\beta_p=1\in B$.
\end{example}

We say that a pair $(X_\bullet, Y_\bullet)$ is  a simplicial pair if $Y_\bullet$ is a simplicial set  and  $X_\bullet$ a simplicial subset of $Y_\bullet$. In other words we have a functor $$(X_\bullet, Y_\bullet): \Delta^{op}\to \Gamma_2^{op}.$$

For a simplicial pair $(X_\bullet, Y_\bullet)$ we define the higher order Hochschild cohomology associated to the triple $(A,B,\varepsilon)$ and a symmetric $A$-bimodule $M$, to be the homology of the complex defined as follows. For every $q\in \mathbb{N}$ we consider $(X_q,Y_q)\in \Gamma_2^{op}$ and take  $C^q_{(X_\bullet, Y_\bullet)}={\mathcal L}((A,B,\varepsilon);M)((X_q,Y_q))$.  We construct a complex by taking the differential induced by the simplicial structure on $(X_\bullet,Y_\bullet)$. More precisely if $d_i:Y_{q+1}\to Y_q$ then we define 
$$\delta^i={\mathcal L}((A,B,\varepsilon);M)(d_i):C^{q}_{(X_\bullet, Y_\bullet)}\to C^{q+1}_{(X_\bullet, Y_\bullet)}$$
and take $\partial_{(X_\bullet, Y_\bullet)}: C^q_{(X_\bullet, Y_\bullet)}\to C^{q+1}_{(X_\bullet, Y_\bullet)}$, 
\begin{equation}
\partial_{(X_\bullet, Y_\bullet)}=\sum_{i=0}^{q+1}(-1)^i\delta^i. \label{eq2}
\end{equation}

\begin{definition} The homology of the above complex is called the higher order Hochschild cohomology associated to the simplicial pair $(X_\bullet, Y_\bullet)$, of the triple $(A,B,\varepsilon)$ with coefficients in $M$ and is denoted by $H^{q}_{(X_\bullet, Y_\bullet)}((A,B,\varepsilon); M)$.
\end{definition}

\begin{remark}
In the event that  $X_\bullet = Y_\bullet$ this definition agrees with the definition of higher order Hochschild cohomology $H^{q}_{X_\bullet}(A, M)$.
\end{remark}

\section{Secondary Cohomology as a Higher Order Cohomology}

In this section we show that when $A$ is commutative and $M$ is a symmetric $A$-bimodule, then the secondary Hochschild cohomology $H^n((A,B,\varepsilon); M)$ is a particular case of the construction from the previous section.

Consider the simplicial pair $(X_\bullet,Y_\bullet)=(S^1,D^2)$, where the sphere $S^1$ is a obtained from the interval $I=[01]$ by identifying the ends of the interval, and the disk $D^2$ is obtained from the $2$-simplex $\Delta=[012]$ by collapsing the edges $[01]$ and $[12]$ (i.e. the boundary of $D^2$ is the edge $[02]$). 

More precisely, we take $X_\bullet$ to be the simplicial set where the only nondegenerate $1$-simplex is $I=[02]$.  We denote by $*_{n}$ the base point in dimension $n$, and by $I^a_b$ the simplex in dimension $n=a+b+1$, where we iterate the $[0]$ vertex $a$ times, and the  $[2]$ vertex $b$ times.  For example, $I^0_0$ is the interval $I$  with $d_0(I^0_0)=d_1(I^0_0)=*_0$, and $I^1_0$ is a $2$-simplex $[002]$ such that $d_0(I^1_0)=d_1(I^1_0)=I^0_0$ and $d_2(I^1_0)=*_1$. 

For $Y_\bullet$, besides the above simplices, we  also have a nondegenerate $2$-simplex  $\Delta=[012]$. Denote it  by $^0\Delta^0_0$ and take  $d_0(^0\Delta^0_0)=d_2(^0\Delta^0_0)=*_1$ and $d_1(^0\Delta^0_0)=I^0_0$.  More generally,  take $^a\Delta^b_c$ the $a+b+c+2$-dimensional simplex obtained by iterating the $[0]$ vertex $a$ times, the $[1]$ vertex $b$ times, and the $[2]$ vertex $c$ times. For example $^1\Delta^0_0$ is a $3$-simplex $[0012]$ with $d_0(^1\Delta^0_0)=d_1(^1\Delta_0^0)= {^0\Delta_0^0}$, $d_2(^1\Delta^0_0)=I_0^1$, and $d_3(^1\Delta^0_0)=\ast_2$. 

In general we have $X^n=\{ *_{n}\}\cup \{ I^a_b|\,a,b\in \mathbb{N},\; a+b=n-1\}$ and $Y^n=X^n\cup \{^a\Delta^b_c|\,a,b,c\in \mathbb{N},\; a+b+c=n-2\}$. The $d_i:Y^n\to Y^{n-1}$ are defined as follows:
\begin{eqnarray}
d_i(*_{n})=*_{n-1},\label{f1}
\end{eqnarray} 
\begin{eqnarray}
d_i(I^a_b)= \left\{\begin{array}{ll}
  *_{a+b}& \mbox{ if $a=0$ and $i=0$}\\ 
  I^{a-1}_b & \mbox{ if  $a\neq 0$ and $i\leq a$ }\\
 I^a_{b-1}& \mbox{ if $b\neq 0$ and $i>a$}\\
 *_{a+b}& \mbox{ if $b=0$ and $i=n=a+1$,}
 \end{array}\right.\label{f2}
\end{eqnarray}
\begin{eqnarray}
d_i(\,^a\Delta^b_c)= \left\{\begin{array}{ll}
  *_{a+b+c+1}& \mbox{ if $a=0$ and $i=0$}\\ 
  ^{a-1}\Delta^b_c & \mbox{ if  $a\neq 0$ and $i\leq a$ }\\
 I^a_{c}& \mbox{ if $b=0$ and $i=a+1$}\\
 ^{a}\Delta^{b-1}_c & \mbox{ if  $b\neq 0$ and $a<i\leq a+b+1$ }\\
 *_{a+b+c+1}& \mbox{ if $c=0$ and $i=n=a+b+2$}\\
 ^{a}\Delta^b_{c-1} & \mbox{ if  $c\neq 0$ and $i\geq a+b+2$}.\\
 \end{array}\right.\label{f3}
 \end{eqnarray}
 
Notice that $I^a_b$ is degenerate if $a+b>0$, and $^a\Delta^b_c$ is degenerate if $a+b+c>0$. 
Also we have that $|X^n|=1+n$ and $|Y^n|=1+n+\frac{n(n-1)}{2}$. In particular we get that 
$${\mathcal L}((A,B,\varepsilon);M)((X^n,Y^n))=Hom_k(A^{\otimes n}\otimes B^{\otimes \frac{n(n-1)}{2}}, M).$$
Next we need to make the identification with the notation from \cite{MR3465889}. First recall that an element in $A^{\otimes n}\otimes B^{\otimes \frac{n(n-1)}{2}}$ was represented by a tensor matrix 

$$
T={\otimes}\left(
\begin{array}{cccccccc}
 a_{1,1}       & \alpha_{1,2} &...&\alpha_{1,n-1}&\alpha_{1,n}\\
 1  & a_{2,2} &...      &\alpha_{2,n-1}&\alpha_{2,n}\\
 .       &. &...&.&.\\
 1& 1  &...&a_{n-1,n-1}&\alpha_{n-1,n}\\
 1& 1  &...&1&a_{n,n}\\
\end{array}
\right),
$$
where $a_{i,i}\in A$ and $\alpha_{i,j}\in B$. 

For $a,b\in \mathbb{N}$ with $a+b+1=n$ the element $I^a_b\in Y^n$ corresponds to the position $(a+1,a+1)$ in the tensor matrix. For $a,b,c\in \mathbb{N}$ with $a+b+c+2=n$ the element $^a\Delta^b_c\in Y^n$ corresponds to the position $(a+1,a+b+2)=(a+1,n-c)$ in the tensor matrix. We also add the symbol $(0,0)$ to correspond  to $*_n$. 

With the above identifications the formulas  (\ref{f2}) and (\ref{f3}) become:
 \begin{eqnarray*}
 d_i((a+1,a+1))= \left\{\begin{array}{ll}
  (0,0)& \mbox{ if $a=0$ and $i=0$}\\ 
  (a,a) & \mbox{ if  $a\neq 0$ and $i\leq a$ }\\
 (a+1,a+1)& \mbox{ if $b\neq 0$ and $i>a$}\\
 (0,0)& \mbox{ if $b=0$ and $i=n=a+1$,}
 \end{array}\right.
 \end{eqnarray*}

\begin{eqnarray*}
d_i(a+1,a+b+2)= \left\{\begin{array}{ll}
  (0,0)& \mbox{ if $a=0$ and $i=0$}\\ 
 (a,a+b+1) & \mbox{ if  $a\neq 0$ and $i\leq a$ }\\
 (a+1,a+1)& \mbox{ if $b=0$ and $i=a+1$}\\
 (a+1,a+b+1)& \mbox{ if  $b\neq 0$ and $a<i\leq a+b+1$ }\\
 (0,0)& \mbox{ if $c=0$ and $i=n=a+b+2$}\\
 (a+1,a+b+2)& \mbox{ if  $c\neq 0$ and $i\geq a+b+2$}.\\
 \end{array}\right. 
 \end{eqnarray*}
If we use the above identification the formula for $\partial_{(S^1,D^2)}$ from  equation (\ref{eq2}) is the same as the formula for differential $\delta^{\varepsilon}_{n-1}$ from equation (\ref{delta}). To summarize we have the following result.

\begin{theorem} Let  $A$ and $B$ be  commutative $k$-algebras, $\varepsilon:B\to A$ a morphism of $k$-algebras  and $M$  a symmetric  $A$-bimodule, then we have
$$ H^{q}((A,B,\varepsilon); M)\simeq H^{q}_{(S^1, D^2)}((A,B,\varepsilon); M). $$
\end{theorem}

\section{Some Remarks}
One can see that  $H^{q}_{(X_\bullet, Y_\bullet)}((A,B,\varepsilon); M)$ is functorial with respect to all of its entries. More precisely, let $(A_1,B_1,\varepsilon_1)$ and $(A_2,B_2, \varepsilon_2)$ be two triples, $M$ a symmetric  $A_2$-bimodule, and $f:A_1\to A_2$ a morphism of $k$-algebras such that $f(B_1)\subseteq B_2$ and  $f\varepsilon_1(b)=\varepsilon_2(f(b))$, then we have the natural morphism $$f^*: H^{q}_{(X_\bullet, Y_\bullet)}((A_2,B_2,\varepsilon_2); M) \to H^{q}_{(X_\bullet, Y_\bullet)}((A_1,B_1,\varepsilon_1); M),$$ 
where the $A_1$-bimodule structure on $M$ is induced by $f$. Also, if $g:M\to N$ is a morphism of symmetric $A$-bimodules then 
 $$g_*: H^{q}_{(X_\bullet, Y_\bullet)}((A,B,\varepsilon); M) \to H^{q}_{(X_\bullet, Y_\bullet)}((A,B,\varepsilon); N).$$ 

Moreover if $h:(X_\bullet, Y_\bullet)\to (Z_\bullet, T_\bullet)$ is a morphism of simplicial sets then 
$$h^*: H^{q}_{(Z_\bullet, T_\bullet)}((A,B,\varepsilon); M) \to H^{q}_{(X_\bullet, Y_\bullet)}((A,B,\varepsilon); M).$$
If we take the natural inclusion of simplicial pairs $i:(S^1,S^1)\to (S^1, D^2)$ with the simplicial structure discussed in the previous section, then $$i^n: H^{n}_{(S^1, D^2)}((A,B,\varepsilon); M)\to H^{n}_{(S^1, S^1)}((A,B,\varepsilon); M)$$ is nothing else but the morphism $\Phi_n:H^n((A,B,\varepsilon);M)\to H^n(A,M)$ discussed in \cite{MR3338544}.

One natural question is whether the construction in this paper depends only of the homotopy type of the geometric realization of the simplicial pair $(X_{\bullet}, Y_{\bullet})$ (or maybe invariant under a certain equivalence relation among simplicial pairs). We explored this problem but we were not able to prove any interesting result. The main issue is finding an equivalence relation among simplicial pairs that is manageable at the  algebraic level.

\section*{Acknowledgment}
The authors would like to thank Andrew Salch for conversations and suggestions about this research.  Bruce would also like to thank his wife Kendall for her continued support.


\addcontentsline{toc}{section}{Bibliography}
\gdef\leftmark{Bibliography} \gdef\rightmark{Bibliography}

\end{document}